\documentclass{article}

\usepackage{times}
\usepackage{graphicx} 
\usepackage{subfigure} 

\usepackage{natbib}

\usepackage{algorithm}
\usepackage{algorithmic}

\usepackage{amsmath, amsthm, amssymb, multirow, paralist}
\newtheorem{thm}{Theorem}
\newtheorem{prop}{Proposition}
\newtheorem{lemma}{Lemma}
\newtheorem{cor}[thm]{Corollary}

 \newtheorem{assumption}{Assumption}
\usepackage{enumitem}
\usepackage{booktabs}       

\usepackage{hyperref}

\usepackage[accepted]{icml2017}

\def \y {\mathbf{y}}

\def \x {\mathbf{x}}

\def \z {\mathbf{z}}
\def \u {\mathbf{u}}

\def \R {\mathbb{R}}
\def \S {\mathcal{S}}

\def \A {\mathcal{A}}

\def \v {\mathbf{v}}

\def \c {\mathbf{c}}

\def \xt {\widetilde{\x}}

\def \a {\mathbf{a}}
\def \b {\mathbf{b}}

\def \xh {\widehat{\x}}

\icmltitlerunning{Convex Constrained  Optimization with Reduced Projections and Improved  Rates}

\begin{document} 

\twocolumn[
\icmltitle{A Richer Theory of Convex Constrained  Optimization\\ with Reduced Projections and Improved  Rates}

\icmlsetsymbol{equal}{*}

\begin{icmlauthorlist}
\icmlauthor{Tianbao Yang}{to}
\icmlauthor{Qihang Lin}{to}
\icmlauthor{Lijun Zhang}{ed}
\end{icmlauthorlist}

\icmlaffiliation{to}{The University of  Iowa, Iowa City, IA 52242, USA}
\icmlaffiliation{ed}{National Key Laboratory for Novel Software Technology, Nanjing University, Nanjing 210023, China}

\icmlcorrespondingauthor{Tianbao Yang}{tianbao-yang@uiowa.edu}

\icmlkeywords{boring formatting information, machine learning, ICML}

\vskip 0.3in
]

 \printAffiliationsAndNotice{}

\begin{abstract} 
This paper focuses on  convex constrained optimization problems, where the solution is subject to a convex inequality constraint. In particular, we aim at challenging problems for which both projection into the constrained domain and a linear optimization under the inequality constraint are  time-consuming, which render both projected gradient methods and conditional gradient methods (a.k.a. the Frank-Wolfe algorithm) expensive. In this paper, we develop projection reduced optimization algorithms for both smooth and non-smooth optimization with improved convergence rates under a certain regularity condition of the constraint function. We first present a general theory of optimization with only one projection. Its application to smooth optimization with only one projection yields $O(1/\epsilon)$ iteration complexity, which  improves over the $O(1/\epsilon^2)$ iteration complexity established before for non-smooth optimization  and can be further reduced under strong convexity. Then we introduce a local error bound condition and develop faster algorithms for non-strongly convex optimization at the price of a logarithmic number of projections. In particular, we achieve  an iteration complexity  of $\widetilde O(1/\epsilon^{2(1-\theta)})$ for non-smooth optimization and $\widetilde O(1/\epsilon^{1-\theta})$ for smooth optimization, where $\theta\in(0,1]$ appearing the local error bound condition characterizes the functional local growth rate around the optimal solutions. Novel applications in solving the constrained $\ell_1$ minimization problem  and a positive semi-definite constrained distance metric learning problem demonstrate that the proposed algorithms achieve significant speed-up compared with previous algorithms. 
\end{abstract} 

\section{Introduction}
 In this paper, we aim at solving the following convex constrained optimization problem: 
\begin{equation}\label{eqn:prob}
\begin{aligned}
\min_{\x\in\R^d}&\quad f(\x),\quad s.t.\;\; c(\x)\leq 0,
\end{aligned}
\vspace*{-0.1in}
\end{equation}
where $f(\x)$ is a smooth or non-smooth convex function and $c(\x)$ is a lower-semicontinuous and convex function. The problem can find applications in machine learning, signal processing, statistics,  marketing optimization, and etc. For example, in distance metric learning one needs to learn a positive semi-definite (PSD) matrix such that similar examples are close to each other and dissimilar examples are far from each other~\citep{Weinberger2006,Xing2003}, where the positive semi-definite constraint can be cast into a convex inequality constraint. Another example arising in compressive sensing is to minimize the $\ell_1$ norm of high-dimensional vector subject to a measurement constraint~\citep{candes-2008-introduction}. Although general interior-point methods can be applied to solve the problem with linear convergence, they suffer from exceedingly high computational cost per-iteration. Another solution is to employ the projected gradient (PG) method~\citep{nesterov2004introductory} or the conditional gradient (CG) method~\citep{frankwolf56}, where the PG method  needs to compute the projection into the constrained domain at each iteration and CG needs to solve a linear optimization problem under the constraint. However, for many constraints (e.g.,  PSD,  quadratic constraints) both projection into the constrained domain and the linear optimization under the constraint are time-consuming, which restrict their capabilities to solving these problems. 

Recently, there emerges a new direction towards addressing the challenge of expensive projection that is to reduce the number of projections. In the seminal paper~\citep{mahdavi-2012-stochstic}, the authors have proposed two algorithms with only one projection at the end of iterations for non-smooth convex  and strongly convex optimization, respectively. The idea of both algorithms is to move the constraint function into the objective function and to control the violation of constraint for intermediate solutions. While their developed  algorithms enjoy an optimal convergence rate for non-smooth optimization (i.e., $O(1/\epsilon^2)$ iteration complexity) and a close-to-optimal convergence rate for strongly convex optimization (i.e., $\widetilde O(1/\epsilon)$~\footnote{where $\widetilde O()$ suppresses a logarithmic factor. }), {\it there still lack of theory and algorithms with reduced projections and faster rates for smooth convex optimization and for convex optimization without strong convexity assumptions}. 

In this paper, we make significant contributions by developing a richer theory of convex constrained optimization with reduced projections and faster rates. To be specific, 
\begin{itemize}[leftmargin=*]
\vspace*{-0.1in}
\item we develop a general framework and theory of optimization with only one projection, where any favorable smooth or non-smooth convex optimization algorithms  can be employed to solve the intermediate augmented unconstrained  objective function.  We discuss in full details the applicability of the proposed algorithms to problems with polyhedral, quadratic or PSD constraints.
\vspace*{-0.1in}

\item Applying the general theory to smooth convex optimization~\footnote{where the constraint function is assumed to be smooth.} with Nesterov's accelerated gradient methods  yields  an iteration complexity of $O(1/\epsilon)$ with only one projection. In addition, when equipped with an optimal algorithm for strongly convex optimization the general theory implies the {\it optimal} iteration complexity of $O(1/\epsilon)$ for strongly convex optimization with only one projection. 
For smooth and strongly convex optimization, the general theory implies an iteration complexity of  $O(1/\epsilon^{\beta})$  where $\beta\in(1/2, 1)$ with only one projection and a sufficiently large number of iterations. 
\vspace*{-0.1in}

\item Building on the general framework and theory, we further develop an improved  theory with faster convergence rates  for non-strongly convex optimization at the price of a logarithmic number of projections. In particular, we show that under a mild local error bound condition, the iteration complexities can be reduced  to $\widetilde O(1/\epsilon^{2(1-\theta)})$ for non-smooth optimization and $\widetilde O(1/\epsilon^{1-\theta})$ for smooth optimization, where $\theta\in(0,1]$ is a constant in the local error bound condition that characterizes the local growth rate of functional values. To our knowledge, these are the best convergence results with only a logarithmic number of projections for non-strongly convex optimization. We also demonstrate their  effectiveness for solving compressive sensing and distance metric learning problems. 

\end{itemize}

\section{Related Work}\label{sec:RW}
The issue of high projection cost in projected gradient descent has received increasing attention in recent years. Most studies are based on the Frank-Wolfe technique that eschews the projection in favor of a linear optimization over the constrained domain~\citep{DBLP:conf/icml/Jaggi13,DBLP:conf/icml/HazanK12,LacosteJulien2013,DBLP:conf/icml/GarberH15}. It happens that for many bounded domains (e.g., bounded balls for vectors and matrices, a PSD constraint with a bounded trace norm)  the linear optimization over the constrained domain is much cheaper than projection into the constrained domain~\citep{DBLP:conf/icml/Jaggi13}. However, there still exist many constraints that render both projection into the constrained domain and linear optimization under the constraint are comparably expensive. Examples include polyhedral constraints, quadratic constraints and a PSD constraint~\footnote{Indeed, a linear optimization over a PSD constraint is ill-posed because the PSD domain is unbounded.}. 

To tackle these complex constraints, the idea of optimization with a reduced number of projections was explored in several studies since~\cite{mahdavi-2012-stochstic}. 
 In a recent paper~\cite{DBLP:conf/uai/ChenYLZC16}, the authors show that for stochastic strongly convex optimization, the optimal convergence rate can be achieved using a logarithmic number of projections. In contrast, the developed theory in this paper implies that only one projection is sufficient to achieve the optimal convergence rate for strongly convex optimization, and a logarithmic number of projections can be used to accelerate convergence rates for non-strongly convex optimization. \citet{DBLP:conf/colt/CotterGP16} proposed a stochastic algorithm for solving heavily constrained problems with many constraint functions by extending the work of~\citep{mahdavi-2012-stochstic}. Nonetheless, their focus is not to improve the convergence rates. \citet{DBLP:conf/icml/ZhangYJH13} studied the smooth and strongly convex optimization and they proposed a stochastic algorithm with $O(\kappa\log(T))$ projections and proved an $O(1/T)$ convergence rate, where $\kappa$ is the condition number and $T$ is the total number of iterations. Nonetheless, if the condition number is high the number of projections could be  very large.  In addition,  their algorithm utilizes the mini-batch to avoid frequent projections in stochastic optimization, which is different from the present paper.


We note that several recent works also exploit different forms of error bound conditions to improve the convergence~\citep{DBLP:journals/jmlr/WangL14,DBLP:journals/corr/So13, DBLP:conf/nips/HouZSL13, DBLP:conf/icml/ZhouZS15,Yang15rsg,DBLP:journals/corr/abs-1607-03815}. Most notably, the technique used in our work is closely related to~\citep{Yang15rsg}. However, for constrained optimization problems the methods in~\citep{Yang15rsg} still need to conduct projections at each iteration. 

Finally, we comment on the differences between the proposed methods and the classical penalty methods that also move the constraint into the objective using a penalty function~\cite{1886529043}. The major differences are that (i) the classical penalty methods typically require solving each subproblem exactly while our methods do not require that; and (ii) the classical penalty methods typically guarantee asymptotic convergence while our methods have explicit convergence rates. 

\section{Preliminaries}\label{sec:Pre}
Let $\Omega = \{\x\in\R^d: c(\x)\leq 0\}$ denote the constrained domain, $\Omega_*$ denote the optimal solution set and $f_*$ denote the optimal objective value.  We denote by $\nabla f(\x)$ the gradient and  by $\partial f(\x)$ the subgradient of a smooth or non-smooth function, respectively. When $f(\x)$ is a non-smooth function, we consider the problem as non-smooth constrained optimization. When both $f(\x)$ and $c(\x)$ are smooth, we consider the problem as smooth constrained optimization. A function $f(\x)$ is $L$-smooth if it has a Lipschitz  continuous gradient, i.e.,  $\|\nabla f(\x) - \nabla f(\y)\|\leq L\|\x - \y\|$, 
where $\|\cdot\|$ denotes the Euclidean norm. A function $f(\x)$ is $\mu$-strongly convex if it satisfies 
$f(\x)\geq f(\y) + \partial f(\y)^{\top}(\x - \y) + \frac{\mu}{2}\|\x - \y\|^2$.

In the sequel, $dist(\x, \Omega)$ denotes the distance of $\x$ to a set $\Omega$, i.e., $dist(\x, \Omega)=\min_{\u\in\Omega}\|\x - \u\|$. Let  $[s]_+$ be a hinge operator that is defined as $[s]_+=s$ if $s\geq 0$, and $[s]_+=0$ if $s<0$. 

Throughout  the paper, we make the  the following assumptions  to facilitate the development of our algorithms and theory.
\begin{assumption}\label{ass:basic} For a convex minimization problem~(\ref{eqn:prob}), we assume
(i) there exists a positive value $\rho>0$ such that 
\begin{align}\label{eqn:keyi0}
\min_{c(\x)=0\atop \v\in \partial c(\x), \v\neq 0}\|\v\|\geq \rho,
\end{align} or more generally  there exists a constant $\rho>0$
 for any $\x\in\R^d$, such that $\x^\natural = \arg\min_{\u\in\R^d, c(\u)\leq 0}\|\u- \x\|^2$ satisfies 
\vspace*{-0.1in} \begin{align}\label{eqn:keyi}
\|\x^\natural - \x\|\leq  [c(\x )]_+/\rho.
\end{align}
(ii) there exists a strictly feasible solution such that $c(\x)<0$;  (iii)  both $f(\x)$ and $c(\x)$  are defined everywhere and are Lipschitz continuous with their Lipschitz constants  denoted by $G$ and $G_c$, respectively. 
\end{assumption}
\vspace*{-0.1in}
We make several remarks about the assumptions. The inequality in~(\ref{eqn:keyi0}) is introduced  in~\citep{mahdavi-2012-stochstic}, which is to ensure the distance from the final solution before projection to constrained domain $\Omega$ is not too large. Note that the inequality in~(\ref{eqn:keyi}) is a more general condition than~(\ref{eqn:keyi0}) as seen from the following lemma. 
\begin{lemma}\label{lem:0}
For any $\x\in\R^d$, let $\x^\natural = \arg\min_{c(\u)\leq 0}\|\u-\x\|^2$. If~(\ref{eqn:keyi0}) holds, then~(\ref{eqn:keyi}) holds. 
\end{lemma}
The above lemma is implicit in the proof of \citep{mahdavi-2012-stochstic}.
We will provide  more discussions about Assumption 1(i) - the key assumption, and exhibit the value of $\rho$ for a number of commonly seen constraints (e.g., polyhedral, quadratic and PSD constraints). To make the presentation more fluent, we postpone these discussions to Section~\ref{sec:dis}.   The strict feasibility assumption (ii) allows us to explore the KKT condition of the projection problem shown below. Assumption (iii) imposes mild Lipschitz continuity conditions on both $f(\x)$ and $c(\x)$. 

Traditional projected gradient descent methods need to solve  the following projection at each iteration 
 $\Pi_{\Omega}[\x] = \arg\min_{ c(\u)\leq 0}\|\u - \x\|^2$. 
Conditional gradient methods (a.k.a. the Frank-Wolfe technique) need to solve the following linear optimization at each iteration
$\min_{\u\in\R^d, c(\u)\leq 0}\u^{\top}\nabla f(\x)$. 
For many constraint functions (see  Section~\ref{sec:dis}), solving the projection problem and the linear optimization could be very expensive.

\section{A General Theory of Optimization with only one projection}\label{sec:Opt}
In this section, we extend the idea of only one projection proposed in~\citep{mahdavi-2012-stochstic} to a general theory, and then present optimization algorithms with only one projection for non-smooth and smooth optimization, respectively. To  tackle the constraint, we introduce a penalty function $h_\gamma(\x)$ parameterized by $\gamma$, which obeys the following certificate: there exist constants $C\geq 0$ and $\lambda>G/\rho$ such that 
\begin{equation}\label{eqn:pen}
\begin{aligned}
&h_\gamma(\x)\geq \lambda [c(\x)]_+,\forall \x\\
&h_\gamma(\x)\leq C\gamma, \forall \x \text{ such that } c(\x)\leq 0.
\end{aligned}
\end{equation}
From the above condition, it is clear that $\gamma\geq 0$. It is notable that the penalty function $h_\gamma(\x)$ will also depend on $\lambda$; however it will be set to a constant value, thus the dependence on $\lambda$ is omitted. We will construct such a penalty function $h_\gamma(\x)$ for non-smooth and smooth optimization in next two subsections. We propose to optimize the following augmented objective function 
\begin{align}\label{eqn:augment}
\min_{\x\in\R^d}F_\gamma(\x) = f(\x) + h_\gamma(\x).
\end{align}
We can employ any applicable optimization algorithms to optimize $F_\gamma(\x)$ pretending that there is no constraint, and finally obtain a solution $\xh_T$ that is not necessarily feasible. In order to obtain a feasible solution, we perform one projection to get $\xt_T = \Pi_{\Omega}(\xh_T)$. 
The following theorem allows us to convert the convergence of $\xh_T$ for $F_\gamma(\x)$ to that of $\xt_T$ for $f(\x)$. 
\begin{thm}\label{thm:1}
Let $\A$ be any iterative optimization algorithm applied to $\min_{\x}F_\gamma(\x)$ with $T$ iterations, which starts with $\x_1$ and returns $\xh_T$ as the final solution, such that the following convergence of $\xh_T$ holds for any $\x\in\R^d$
\begin{align}\label{eqn:key0}
F_\gamma(\xh_T) - F_\gamma(\x)\leq B_T(\gamma; \x, \x_1),
\end{align}
where $B_T(\gamma; \x, \x_1)\rightarrow 0$ when $T\rightarrow \infty$. Suppose that Assumption~\ref{ass:basic} hold,  then
\begin{align}\label{eqn:conv}
f(\xt_T) - f(\x_*) \leq \frac{\lambda\rho}{\lambda\rho - G}(C\gamma+B_T(\gamma; \x_*, \x_1)),
\end{align}
where $\xt_T= \Pi_{\Omega}[\xh_T]$ and $\x_*$ is an optimal solution to~(\ref{eqn:prob}). 
\end{thm}
{\bf Remark:} It is worth mentioning that we omit some constant factors in the convergence bound $B_T(\gamma;  \x, \x_1)$ that are irrelevant to our discussions. The notation $B_T(\gamma; \x, \x_1)$ emphasizes that it is a function of $\gamma$ and depends on $\x_1$ and a target solution $\x$ and it will be referred to as $B_T$.  In the next several subsections, we will see that by carefully choosing the penalty function $h_\gamma(\x)$ we are able to provide nice convergence for smooth and non-smooth optimization with only one projection.  In the above theorem, we assume the optimization algorithm $\mathcal A$ is deterministic. However, a similar result can be easily extended to a stochastic optimization algorithm $\mathcal A$. 
\vspace{-0.1in}\begin{proof}
First, we consider $c(\xh_T)\leq 0$, which implies that $\xh_T= \xt_T$. Due to the certificate of $h_\gamma(\x)$,  $F_\gamma(\widetilde\x_T)\geq f(\widetilde\x_T)$ and $F_\gamma(\x_*) \leq  f(\x_*)  + C\gamma$. Hence $f(\xt_T)\leq F_\gamma(\xh_T)\leq F_\gamma(\x_*)+B_T(\gamma; \x_1, \x_*) \leq f(\x_*) + C\gamma + B_T(\gamma; \x_1, \x_*)$. Then~(\ref{eqn:conv}) follows due to $\lambda\rho/(\lambda \rho - G)\geq 1$. Next, we assume $c(\widehat\x_T)>0$. 
Inequality~(\ref{eqn:key0}) implies that
\begin{equation}\label{eqn:B0}
f(\widehat\x_T) + \lambda [c(\xh_T)]_+\leq f(\x_*) + C\gamma  +  B_T(\gamma; \x_*,  \x_1).
\end{equation}
By Assumption~\ref{ass:basic}(i), we have
$[c(\widehat\x_T)]_+\geq \rho\|\widehat\x_T - \widetilde\x_T\|$. 
Combined  with~(\ref{eqn:B0}) we have
\begin{align*}
\lambda\rho&\|\widehat\x_T - \widetilde\x_T\|\leq f(\x_*) - f(\widehat\x_T)+ C\gamma + B_T(\gamma;  \x_*,\x_1)\\
&\leq G\|\widehat\x_T - \widetilde\x_T\| + C\gamma + B_T(\gamma; \x_*,  \x_1),
\end{align*}
where the last inequality follows that fact $f(\x_*) - f(\widehat\x_T)\leq f(\x_*) - f(\widetilde\x_T) + f(\widetilde\x_T) - f(\widehat\x_T) \leq G\|\widehat\x_T - \widetilde\x_T\|$ because the Lipschitz property and $f(\x_*)\leq f(\xt_T)$.
Therefore we have
\[
\|\widehat\x_T - \widetilde\x_T\| \leq \frac{C\gamma+ B_T(\gamma;  \x_*, \x_1,)}{\lambda\rho - G}.
\]
Finally, we obtain
\begin{align*}
f(\widetilde\x_T) - f(\x_*)&\leq f(\widetilde\x_T) - f(\widehat\x_T)+ f(\widehat\x_T) - f(\x_*)\\
&\leq G \|\widehat\x_T - \widetilde\x_T\| + C\gamma+ B_T(\gamma;  \x_*, \x_1)\\
&\leq \frac{\lambda \rho}{\lambda\rho - G}(C\gamma+ B_T(\gamma;\x_*,  \x_1)).
\end{align*}
\end{proof}

\subsection{Non-smooth Optimization}
\vspace*{-0.05in}
Since an optimal convergence rate for general non-smooth optimization with only one projection has been attained in~\citep{mahdavi-2012-stochstic},  in this subsection we present an optimal convergence  result for strongly convex problems.   For non-smooth optimization, we can choose $$ h(\x)=\lambda [c(\x)]_+,$$ and hence $\gamma=0$. We will use deterministic subgradient descent as an example to demonstrate the convergence for $f(\x)$, though many other optimization algorithms designed for non-smooth optimization are applicable (e.g., the stochastic subgradient method). 
The update of subgradient descent method is given by the following
\begin{align}\label{eqn:sg}
\x_{t+1}  = \x_{t}  - \eta_t \partial F(\x_{t}), \quad  t=1,\ldots, T,
\end{align}
where $\eta_t$ is an appropriate step size. If $f(\x)$ is $\mu$-strongly convex, the step size can be set as $\eta_t=1/(\mu t)$ and the final solution can be computed 
by the $\alpha$-suffix averaging 
$\xh_T = \frac{1}{\alpha T}\sum_{t=(1-\alpha)T+ 1}^T\x_t$ with $\alpha>0$~\citep{RakhlinSS12}, 
or by the polynomial decay averaging with $\xh_t = (1 - \frac{s+1}{s + t})\xh_{t-1} + \frac{s+1}{s + t}\x_t$ and $s\geq 1$~\citep{DBLP:conf/icml/Shamir013}. Both schemes can attain  $B_T=O(1/(\mu T))$ for  the convergence of $F(\x)$ when $f(\x)$ is $\mu$-strongly convex.  Combining this with Theorem~\ref{thm:1}, we have the following convergence result with the proof omitted due to its simplicity. 
\begin{cor}\label{cor:1}
Suppose that Assumption~\ref{ass:basic} holds and $f(\x)$ is $\mu$-strongly convex. Set $F(\x) = f(\x) + \lambda [c(\x)]_+$ with $\lambda\geq G/\rho$. Let~(\ref{eqn:sg}) run for $T$ iterations with 
$\eta_t = 1/(\mu t)$. Let   $\xh_T$ be computed by $\alpha$-suffix averaging or the polynomial decay averaging. Then with only one projection $\xt_T=\Pi_{\Omega}(\xh_T)$, we achieve
\[
f(\xt_T) - f_* \leq \frac{\lambda\rho}{\lambda \rho -G}\frac{(G + \lambda G_c)^2O(1)}{\mu T}.
\]
\end{cor}
{\bf Remark:} 
We note  that the $O(1/(\mu T))$ is also achieved for strongly convex optimization in~\cite{DBLP:conf/icml/ZhangYJH13,DBLP:conf/uai/ChenYLZC16} but with a logarithmic number of projections. In contrast, Corollary~\ref{cor:1} implies only one projection is sufficient to achieve the optimal convergence for strongly convex optimization. 


\subsection{Smooth Optimization}\label{subsec:so}
\vspace*{-0.05in}
For smooth optimization, we consider both $f(\x)$ and $c(\x)$ to be smooth~\footnote{it can be extended to when $f(\x)$ is non-smooth but its proximal mapping can be easily solved.}. Let the smoothness parameter of $f(\x)$ and $c(\x)$ be $L_f$ and $L_c$, respectively. In order to ensure the augmented function $F_\gamma(\x)$ to be still a smooth function, we construct the following penalty function 
\begin{align}\label{eqn:pens}
h_\gamma(\x) = \gamma \ln\left(1 + \exp\left(\lambda c(\x)/\gamma\right)\right).
\end{align}
The following proposition shows that $h_\gamma(\x)$ is a smooth function and obeys the condition in~(\ref{eqn:pen}).
\begin{prop}\label{prop:1}
 Suppose $c(\x)$ is $L_c$-smooth and $G_c$-Lipschitz continuous. The penalty function in~(\ref{eqn:pens}) is a $(\lambda L_c + \frac{\lambda^2G_c^2}{4\gamma})$-smooth function and satisfies (i) $h_\gamma(\x)\geq \lambda [c(\x)]_+$ and (ii)
$h_\gamma(\x)\leq \gamma\ln 2$, $\forall \x$ such that $c(\x)\leq 0$. 
\end{prop}
 Then $F_\gamma(\x)$  is a smooth function and its  smoothness parameter is given by 
$
L_F = L_f + \lambda L_c + \frac{\lambda^2G_c^2}{4\gamma}$. 
Next, we will establish the convergence for $f(\x)$ using  Nesterov's optimal accelerated gradient (NAG) methods. The update of  one variant of NAG can be written as follows
\begin{equation}\label{eqn:NAG}
\begin{aligned}
\x_{t+1} &= \y_t - \nabla F_\gamma(\y_t)/L_F\\
\y_{t+1} & = \x_{t+1} + \beta_{t+1}(\x_{t+1} - \x_t),
\end{aligned}
\end{equation}
where the value of $\beta_t$ can be set to different values depending on whether $f(\x)$ is strongly convex or not (see Corollary~\ref{cor:2}). Previous work have established the convergence of $\xh_T = \x_{T}$ for $F_\gamma(\x)$, in particular $B_T =O( \frac{L_F}{T^2})$ for smooth non-strongly convex optimization and $B_T = O\left(L_F\exp(- T\sqrt{\frac{\mu}{L_F}})\right)$ for smooth and strongly convex optimization. By combining these results with Theorem~\ref{thm:1} and appropriately setting $\gamma$, we can achieve the following convergence of $\xt_T$ for $f(\x)$. 
\begin{cor}\label{cor:2}
Suppose that Assumption~\ref{ass:basic} holds, $dist(\y_0, \Omega_*)\leq D$, $f(\x)$ is $L_f$-smooth and $c(\x)$ is $L_c$-smooth. Set $F_\gamma(\x) = f(\x) + h_\gamma(\x)$ with  $\lambda>G/\rho$ and $h_\gamma(\x)$ being~(\ref{eqn:pens}). 
Let~(\ref{eqn:NAG}) run for $T$ iterations and $\xt_T= \Pi_{\Omega}(\x_T)$. 
\begin{itemize}[leftmargin=*,topsep=0pt,itemsep=-1ex,partopsep=1ex,parsep=1ex, itemindent=2ex]
\item If $f(\x)$ is convex, we can set $\gamma = \frac{\lambda G_cD}{(T+1)\sqrt{2\ln 2}}$, $\beta_t = \frac{\tau_{t-1}-1}{\tau_t}$, where $\tau_t = \frac{1+\sqrt{1+4\tau_{t-1}^2}}{2}$ with $\tau_0=1$, and achieve 
\vspace*{-0.1in}
\[
f(\xt_T) - f_* \leq \frac{\lambda \rho}{\lambda \rho - G}\left[\frac{\lambda G_cD\sqrt{2\ln 2}}{T+1 } + \frac{(L_f+\lambda L_c) D^2}{(T+1)^2}\right]
\]  
\item If $f(\x)$ is $\mu$-strongly convex, we can set $\gamma = \frac{1}{T^{2\alpha}}$ with $\alpha\in(1/2, 1)$ and $\beta_t = \frac{\sqrt{L_F} - \sqrt{\mu}}{\sqrt{L_F} + \sqrt{\mu}}$, and achieve 
\[
f(\xt_T) - f_* \leq O\left(\frac{1}{T^{2\alpha}} + \frac{1}{T^{4\alpha}}\right),
\]  
\vspace*{-0.1in}
as long as $T\geq \left(\frac{L_f + \lambda L_c + \lambda^2 G_c^2/4}{\mu}\right)^{\frac{1}{2(1-\alpha)}}(4\alpha\ln T)^{\frac{1}{1-\alpha}}$. 
\end{itemize}
\end{cor}
{\bf Remark:} The convergence results above indicate an  $O(1/\epsilon)$ iteration complexity for smooth optimization and $O(1/\epsilon^{1/(2\alpha)})$ with $\alpha\in(1/2, 1)$ for smooth and strongly convex optimization with only one projection. Omitted proofs can be found in appendix.

\section{Improved Convergence for Non-strongly Convex Optimization}\label{sec:Imp}
In this section, we will develop improved convergence for non-strongly convex optimization at a price of a logarithmic number of projections by considering an additional condition on the target problem.  
To facilitate the presentation, we first introduce some notations.  
The $\epsilon$-sublevel set  $\S_\epsilon$  and $\epsilon$-level set  $\mathcal L_\epsilon$ of the problem~(\ref{eqn:prob}) are denoted by 
$\S_\epsilon = \{\x\in\Omega: f(\x) \leq f_* + \epsilon\}$, and $\mathcal L_\epsilon = \{\x\in\Omega: f(\x) = f_* + \epsilon\}$, respectively. 
Let $\x_{\epsilon}^\dagger$ denote the closest point in the $\epsilon$-sublevel set $\S_\epsilon$ to $\x\in\Omega$, i.e.,
\begin{equation}\label{eqn:ec}
\begin{aligned}
\x_\epsilon^{\dagger}&=\arg\min_{\u\in\Omega}\|\u - \x\|^2,\quad \text{s.t.}\quad f(\u)\leq f_* + \epsilon.
\end{aligned}
\end{equation}
Let $\x^*$ denote the closest optimal solution in $\Omega_*$ to $\x$, i.e., $\x^* = \arg\min_{\u\in\Omega_*}\|\u - \x\|^2$. 

In this section, we will make the following additional assumption about the problem~(\ref{eqn:prob}). 
\begin{assumption}\label{ass:rsg} For a convex minimization problem~(\ref{eqn:prob}), we assume
(i) there exist $\x_0\in\Omega$ and $\epsilon_0\geq 0$ such that $f(\x_0) - \min_{\x\in\Omega}f(\x)\leq \epsilon_0$;
(ii) $\Omega_*$ is a non-empty convex compact set; 
(iii) the optimization problem~(\ref{eqn:prob}) satisfies a local error bound condition, i.e., there exist $\theta\in(0,1]$ and $\sigma>0$ such that for any $\x\in\S_\epsilon$ we have $dist(\x, \Omega_*)\leq \sigma (f(\x) - f_*)^{\theta}$ where $\Omega_*$ denotes the optimal set and $f_*$ denotes the optimal value. 
\end{assumption}
\vspace*{-0.10in}{\bf Remark:} we would like to remark that the new assumption only imposes mild conditions on the problem. In particular,  Assumption~\ref{ass:rsg} (i) supposes there is a lower bound of the optimal value $f_*$, which usually holds in machine learning applications where the objective function if non-negative;  Assumption~\ref{ass:rsg} (ii) ensures that $\mathcal S_\epsilon$ is also bounded~\citep{rockafellar1970convex}, therefore the $\sigma$ in the local error bound is finite, which can be easily satisfied for a norm regularized or constraint problems; the local error bound condition holds for a broad family of functions (e.g., semi-algebraic functions or real subanalytic functions~\cite{arxiv:1510.08234,Yang15rsg}).  In Section~\ref{sec:app}, we will  also demonstrate several applications of the improved algorithms proposed in this section by establishing the local error bound condition. 

Although the local error bound condition is much weaker than the strong convexity assumption, below we will propose novel algorithms leveraging this condition with faster convergence and only a logarithmic number of projections.

\subsection{Non-smooth Optimization}
\vspace*{-0.05in}
To establish an improved convergence for non-smooth optimization, we develop a new algorithm shown in Algorithm 1 based on subgradient descent (GD) method, to which we refer as LoPGD. The algorithm runs for $K$ epochs and each epoch employs GD for minimizing $F(\x)=f(\x)+\lambda[c(\x)]_+$ with a feasible solution $\x_{k-1}\in\Omega$ as a starting point and $t$ iterations of updates. At the end of each epoch, the averaged solution $\xh_k$ is projected into the constrained domain $\Omega$ and the solution $\x_k$ will be used as the starting point for next epoch. The step size $\eta_k$ is decreased by half every epoch starting from a given value $\eta_1$. The theorem below establishes the iteration complexity of LoPGD  and also exhibits the values of $K$, $t$ and $\eta_1$. To simplify notations, we let $p = \frac{\lambda\rho}{\lambda \rho - G}$ and $\bar G = G + \lambda G_c$. 
\begin{thm}\label{thm:2}
Suppose Assumptions~\ref{ass:basic} and~\ref{ass:rsg} hold. Let $\eta_1 = \frac{\epsilon_0}{2p \bar G^2}$, $K = \lceil \log_2(\epsilon_0/\epsilon)\rceil$ and $t= \frac{4\sigma^2p^2\bar G^2}{\epsilon^{2(1-\theta)}}$ in Algorithm~\ref{fig:1}, where $\theta$ and $\sigma$ are constants appearing in the local error bound condition. Then  $f(\x_K) - f_* \leq 2\epsilon$. 
\end{thm}
\vspace*{-0.1in}
{\bf Remark:} Since the projection is only conducted at the end of each epoch and the total number of epochs is at most $K =\lceil \log_2(\epsilon_0/\epsilon)\rceil$, so the total number of projections is only a logarithmic number $K$. The iteration complexity in~Theorem~\ref{thm:2} is $\widetilde O(1/\epsilon^{2(1-\theta)})$ that improves the standard result of $O(1/\epsilon^2)$ without strong convexity. With $\theta=1/2$, we can achieve $\widetilde O(1/\epsilon)$ iteration complexity with only $O(\log(1/\epsilon))$ projections. 
 
 \setlength\floatsep{0.2\baselineskip plus 3pt minus 2pt}
\setlength\textfloatsep{0.2\baselineskip plus 3pt minus 2pt}
\setlength\intextsep{0.2\baselineskip plus 3pt minus 2 pt}
\begin{algorithm}[t]
\caption{LoPGD} \label{fig:1}
\begin{algorithmic}[1]
\STATE \textbf{INPUT}:  $K\in\mathbb N^+$ , $t\in\mathbb N^+$, $\eta_1$
\STATE \textbf{Initialization}: $\x_0\in \Omega, \epsilon_0$
 \FOR{$k=1, 2, \ldots, K$}
   \STATE Let $\x^k_1=\x_{k-1}$
   \FOR{$s=1, 2, \ldots, t-1$}
        \STATE Update $\x^{k}_{s+1} = \x^k_{s} - \eta_k \partial F(\x^k_{s})$
   \ENDFOR
   
   \STATE Let $\xh_k = \sum_{s=1}^t\x^k_{s}/t$
   \STATE Let $\x_k = \Pi_{\Omega}[\xh_k]$ and $\eta_{k+1}=\eta_k/2$
\ENDFOR
\end{algorithmic}
\end{algorithm}
\begin{algorithm}[t]
\caption{LoPNAG} \label{fig:1-1}
\begin{algorithmic}[1]
\STATE \textbf{INPUT}:  $K\in\mathbb N^+$ , $t_1,\ldots, t_K\in\mathbb N^+$, $\gamma_1$
\STATE \textbf{Initialization}: $\x_0\in \Omega, \epsilon_0$
 \FOR{$k=1, 2, \ldots, K$}
   \STATE Let $\y^k_0=\x_{k-1}$
   \FOR{$s=0, 1, 2, \ldots, t_k-1$}
        \STATE Update $\x^{k}_{s+1} = \y^k_{s} - \frac{1}{L_{k}} \nabla F_{\gamma_k}(\x^k_{s})$
        \STATE Update $\y^k_{s+1}=  \x^k_{s+1} + \beta_{s+1}(\x^k_{s+1} - \x_s^k)$
   \ENDFOR
   
   \STATE Let $\xh_k = \x^k_{t_k}$, $\x_k = \Pi_{\Omega}[\xh_k]$ and  $\gamma_{k+1}=\gamma_k/2$

\ENDFOR
\end{algorithmic}
\end{algorithm}

\subsection{Smooth Optimization}
\vspace*{-0.05in}
Similar to non-smooth optimization, we also develop a new algorithm based on NAG shown in Algorithm 2, where $F_{\gamma}(\x)$ is defined using $h_{\gamma}(\x)$ in~(\ref{eqn:pens}), $L_k = L_{F_{\gamma_k}}$ is the smoothness parameter of $F_{\gamma_k}$ and $\beta_{s}=\frac{\tau_{s-1} - 1}{\tau_s}, s=1,\ldots, $ is a sequence with $\tau_s$ updated as in Corollary~\ref{cor:2}. We refer to this algorithm as LoPNAG. The key idea is to use to a sequence of reducing values for $\gamma_k$ instead of using a small value as in Corollary~\ref{cor:2}, and solve each augmented unconstrained problem $F_{\gamma_k}(\x)$ approximately with one projection.  The theorem below exhibits the iteration complexity of LoPNAG and reveals the values of $K$, $\gamma_1$ and $t_1, \ldots, t_K$. To simplify notations, we let $\bar L = L_f + \lambda L_c$. 
\begin{thm}\label{thm:3}
Suppose Assumptions~\ref{ass:basic} and~\ref{ass:rsg} hold and $f(\x)$ is $L_f$-smooth and $c(\x)$ is $L_c$-smooth. Let $\gamma_1 = \frac{\epsilon_0}{6p\ln 2}$, $K = \lceil \log_2(\epsilon_0/\epsilon)\rceil$ and $t_k=\frac{\sigma}{\epsilon^{1-\theta}}\max \{\lambda G_cp\sqrt{18\ln 2}, \sqrt{12(L_f + \lambda L_c)\epsilon_{0}/2^{k-1}}\}$ in Algorithm~\ref{fig:1-1}, where $\theta$ and $\sigma$ are constants appearing in the local error bound condition. Then $f(\x_K) - f_* \leq 2\epsilon$. 
\end{thm}\vspace*{-0.1in}
{\bf Remark: } It is not difficult to show that the total number of iterations is bounded by $\widetilde O(1/\epsilon^{1-\theta})$, which improves the one in Corollary~\ref{cor:2} without strong convexity.  If $f(\x)$ is a simple non-smooth function whose proximal mapping can be easily computed (e.g., $\ell_1$ norm), we can replace step 6 in Algorithm~\ref{fig:1-1} by a proximal mapping to handle $f(\x)$, which gives the same convergence result in Theorem~\ref{thm:3}. An example is presented in Section~\ref{sec:app} for compressive sensing with $\theta=1/2$.

\section{Discussion of Assumption~\ref{ass:basic} (i)}\label{sec:dis}
One might note that a key condition for developing  the theory with reduced projections is Assumption~\ref{ass:basic} (i). Although \citet{mahdavi-2012-stochstic} has briefly mentioned that the condition can be satisfied for a PSD cone or a Polytope (a bounded polyhedron), their discussion lacks of details in particular on the value of $\rho$ in~(\ref{eqn:keyi0}) or~(\ref{eqn:keyi}). 
Below, we discuss the condition in details about three types of constraints. 

\vspace*{-0.1in}
\paragraph{Polyhedral constraints.}
First, we show that when $c(\x)$ is a polyhedral function, i.e., its epigraph is a polyhedron (not necessarily bounded), the inequality~(\ref{eqn:keyi}) is satisfied. To this end, we explore the polyhedral error bound (PEB) condition~\citep{DBLP:journals/mp/GilpinPS12,Yang15rsg}. In particular, if we consider an optimization problem, 
$\min_{\x\in\R^d}h(\x)$, 
where the epigraph of $h(\x)$ is polyhedron. Let $\mathcal H_*$ denote the optimal set and $h_*$ denote the optimal value of the problem above. The PEB says that there exists $\rho>0$ such that for any $\x\in\R^d$
\begin{align}\label{eqn:PEB}
dist(\x, \mathcal H_*)\leq (h(\x) -h_*)/\rho.
\end{align}
To show that the inequality~(\ref{eqn:keyi}) holds for a polyhedral function $c(\cdot)$, we can consider the  optimization problem  $\min_{\x\in\R^d} [c(\x)]_+$. 
The optimal set of the above problem is given by $\mathcal H_* =\{\x\in\R^d: c(\x)\leq 0\}$. For any $\x$ such that $c(\x)>0$,  let $\x^\natural=  \arg\min_{c(\u)\leq 0}\|\u-\x\|^2 $ be the closest point in the optimal set to $\x$. Therefore  if $c(\cdot)$ is a polyhedral function so does $[c(\x)]_+$,  by the PEB condition~(\ref{eqn:PEB}) there exists a $\rho>0$ such that
$\|\x - \x^\natural\|\leq ([c(\x)]_+  - \min_\x [c(\x)]_+)/\rho = [c(\x)]_+/\rho.$ 
Let us consider a concrete example, where the problem has  a set of affine inequalities   $\c_i^{\top}\x - b_i\leq 0, i =1,\ldots, m$. There are two methods to encode this into a single constraint function $c(\x)\leq 0$. The first method is to use $c(\x) = \max_{1\leq i\leq m}\c_i^{\top}\x - b_i$, which is a polyhedral function and therefore satisfies~(\ref{eqn:keyi}). The second method is to use $c(\x) = \|[C\x - \b]_+\|$, where $[\a]_+=\max(0, \a)$ and $C=(\c_1,\ldots, \c_m)^{\top}$. Thus $[c(\x)]_+ =  \|[C\x - \b]_+\|$. The inequality~(\ref{eqn:keyi}) is then guaranteed by Hoffman's bound and  the parameter $\rho$ is given by the minimum non-zero eigenvalue of $C^{\top}C$~\citep{DBLP:journals/jmlr/WangL14}. Note that the projection onto a polyhedron  is a linear constrained quadratic programming problem, and the linear optimization over a polyhedron is a linear programming problem. Both have polynomial time complexity that would be high if $m$ and $d$ are large~\citep{Karmarkar:1984:NPA:800057.808695,zbMATH03746828}.  


\paragraph{Quadratic constraint.} A quadratic constraint can take the form of $\|A\x - \y\|^2\leq \tau$, where $A\in\R^{m\times d}$ and $\y\in\R^m$. Such a constraint appears in compressive sensing~\citep{candes-2008-introduction}\footnote{Here we use the square constraint to make it a smooth function so that the proposed algorithms for smooth optimization are applicable by using proximal gradient mapping to handle the $\ell_1$ norm.}, where the goal is to reconstruct a sparse high-dimensional vector $\x$ from a small number of noisy measurements $\y = A\x + \varepsilon\in\R^m$ with $m\ll d$. The corresponding optimization problem is 
\begin{equation}\label{eqn:cs}
\begin{aligned}
\textstyle \min_{\x\in\R^d}&\quad\|\x\|_1,\quad s.t. \;\;\|A\x - \y\|^2\leq \tau.  
\end{aligned}
\end{equation}
where $\tau\geq \|\varepsilon\|^2$ is an upper bound on the magnitude of the noise. To check the Assumption~\ref{ass:basic}(i), we note that $c(\x) = \|A\x-\y\|^2 - \tau$ and $\nabla c(\x) = A^{\top}(A\x - \y)$. Let us consider that $A$ has a full row rank~\footnote{which is reasonable because $m\ll d$.} and denote by $\v = A\x - \y$, then on the boundary $c(\x)=0$ we have $\|\v\|=\sqrt{\tau}$ and $\|A^{\top}\v\|\geq \sqrt{\tau\lambda_{\min}(AA^{\top})}$, where $\lambda_{\min}(AA^{\top})>0$ is the minimum eigenvalue of $AA^{\top}\in\R^{m\times m}$. Therefore the Assumption~\ref{ass:basic}(i) is satisfied with $\rho  = \sqrt{\tau\lambda_{\min}(AA^{\top}) }$. It is notable that the projection  and the linear optimization under the quadratic constraint require solving a quadratic programming problem and therefore could be expensive. 
\vspace*{-0.1in}

\paragraph{PSD constraint.} A PSD constraint  $X\succeq0$ for $X\in\R^{d\times d}$ can be written as an inequality constraint $-\lambda_{\min}(X) \leq 0$, where $\lambda_{\min}(X)$ denotes the minimum eigen-value of $X$. The subgradient of $c(X) =  - \lambda_{\min}(X)$ when $\lambda_{\min}(X)=0$ is given by $\text{Conv}\{-\u\u^{\top}| \|\u\|= 1, X\u = 0\}$
, i.e.,  the convex hull of the outer products of normalized vectors in the null space of  the matrix $X$. 
In appendix, we show that if the dimension of the null space of $X$ is $r$ with $1\leq r\leq d$,  the norm of the subgradient of $c(X)$ on the boundary $c(X)=0$ is lower bounded by $\rho = \frac{1}{\sqrt{r}}\geq \frac{1}{\sqrt{d}}$. 
Finally, we note that computing a subgradient of $[c(X)]_+$ only needs to compute one eigen-vector corresponding to the smallest eigen-value.  In contrast,  both projection and linear optimization under a PSD constraint could be very expensive for high-dimensional problems. In particular, the projection onto a PSD domain needs to conduct a singular value decomposition. The linear optimization over a PSD cone is ill-posed due to that PSD cone is not compact (the solution is either 0 or infinity). One may add an artificial constraint on the upper bound of the eigen-values. According to~\cite{DBLP:conf/icml/Jaggi13}, the time complexity for solving this linear optimization problem approximately up to an accuracy level $\epsilon'$ is $O(Nd^{1.5}/{\epsilon'}^{2.5})$ with $N$ being the number of non-zeros in the gradient and $\epsilon'$ decreasing iteratively required in the Frank-Wolfe method,  which could be much more expensive especially for high-dimensional problems and in later iterations than computing the first eigen-pairs at each iteration in our methods. 
\section{Applications}\label{sec:app}
\vspace*{-0.05in}
\subsection{Compressive Sensing}
\vspace*{-0.05in}
\begin{table*}[t]
{\centering\caption{LoPNAG vs. NESTA for solving the compressive sensing problem. 
}}
\label{tab:cs}
\centering{\small
\begin{tabular}{llll|llll}
\toprule
\multicolumn{4}{c|}{LoPNAG}&\multicolumn{4}{c}{NESTA}\\\hline
Iters - Projs &Rec. Err.& Objective& Time (s)  &Iters - Projs&Rec. Err.	& Objective& Time (s) \\ \hline
$5000$ - $1$	& $0.018017$& $52.042878$	& $18.04$&$1000$ - $2000$&$0.137798$&$52.703275$&$48.49$ \\
$10000$ - $2$ & $0.018038$&$52.042418$ & $35.88$&$3000$ - $6000$&$0.018669$&$52.050051$&$93.84$\\
$15000$ - $3$ & $0.018043$&${52.042358}$& ${53.09}$&$5000$ - $10000$&$0.018659$&$52.050046$&$245.23$\\
$20000$ - $4$ & $0.018043$&$52.042358$& $70.24$&$8000$ - $16000$&$0.018657$&$52.050045$&$404.72$\\
$25000$ - $5$ & $0.018043$&$52.042358$& $87.32$&$10000$ - $20000$&$0.018657$&$52.050044$&$501.65$\\
\bottomrule
\end{tabular}}
\vspace*{-0.2in}
\end{table*}
\begin{table*}[t]
\caption{LoPGD vs. OPGD and PGD for solving the considered distance metric learning problem. }
\label{tab:dml}
\centering{\small
\begin{tabular}{lll|lll|lll}
\toprule
\multicolumn{3}{c|}{LoPGD}&\multicolumn{3}{c|}{OPGD}&\multicolumn{3}{c}{PGD}\\\hline
Iters - Projs & Objective& Time (h)  &Iters -  Projs	& Objective& Time (h) &Iters - Projs	& Objective& Time (h) \\ \hline
$1000$ - $1$	& $0.0953$	& $0.22$&$1000$ - $1$&$0.1707$&$0.20$&$1000$ - $1000$&$0.1491$&$7.97$ \\
$2000$ - $2$	& $0.0695$	& $0.43$&$2000$ - $1$&$0.1583$&$0.40$&$2000$ - $2000$&$0.1278$&$15.46$ \\
$4000$ - $4$ &$0.0494$ & $0.87$&$4000$ - $1$&$0.1469$&$0.80$&$4000$ - $4000$&0.1072&29.39\\
$6000$ - $6$ &$0.0428$& $1.33$&$6000$ - $1$&$0.1398$&$1.22$&$6000$ - $6000$&$0.0957$&$43.36$\\
$8000$ - $8$ &$0.0405$& $1.89$&$8000$ - $1$&$0.1343$&$1.64$&$8000$ - $8000$ & $0.0879$&$57.43$\\
\bottomrule
\end{tabular}}
\vspace*{-0.2in}
\end{table*}
We first consider a compressive sensing problem in~(\ref{eqn:cs}). \citet{Becker:2011:NFA:2078698.2078702}  proposed an optimization algorithm  based on the Nesterov's smoothing  and the Nesterov's optimal method for the smoothed problem,  known as NESTA. It needs to perform the projection into the domain $\|A\x - \y\|^2\leq \tau$ at every iteration and has an iteration complexity of $O(1/\epsilon)$. In contrast, the presented algorithm with {\it only one projection} in Section~\ref{subsec:so} using Nesterov's accelerated {\it proximal} gradient method~\cite{Beck:2009:FIS:1658360.1658364} to solve the unconstrained problem enjoys an iteration complexity of $O(1/\epsilon)$. Moreover, we present a theorem below showing that the problem~(\ref{eqn:cs}) satisfies the local error bound condition with $\theta=1/2$, and hence the presented LoPNAG enjoys an $\widetilde O(1/\sqrt{\epsilon})$ iteration complexity with only a logarithmic number of projections. 
\begin{thm}\label{thm:cs}
Let $f(\x) =\|\x\|_1, c(\x) = \|A\x - \y\|^2 - \tau$, $\Omega_*$ denote the optimal set and $f_*$ be the optimal solution to~(\ref{eqn:cs}). Assume that there exists $\x_0$ such that $\|A\x_0-\y\|^2< \tau$ and $0\not\in\Omega_*$. Then for any $\epsilon>0$, $\x\in\R^d$ such that $c(\x)\leq 0$ and $f(\x)\leq  f_* +\epsilon$, there exists $0<\sigma<\infty$ such that $dist(\x, \Omega_*)\leq \sigma (f(\x) - f_*)^{1/2}$. Hence, LoPNAG can have an iteration complexity of $\widetilde O(1/\sqrt{\epsilon})$ with only $O(\log(1/\epsilon))$ projections. 
\end{thm}
\vspace*{-0.1in}
Next, we demonstrate the effectiveness of  the LoPNAG for solving the compressive sensing problem in~(\ref{eqn:cs}) by comparing with NESTA. We generate a synthetic data for testing. In particular, we generate a random measurement matrix $A\in\R^{m\times d}$ with $m = 1000$ and $d = 5000$. The entries of the matrix $A$ are generated independently with the uniform distribution over the interval $[-1,+1]$. The vector $\x_* \in \R^d$ is generated with the same distribution at $100$ randomly chosen coordinates. The noise $\varepsilon\in\R^m$ is a dense vector with independent random entries with the uniform distribution over the interval $[-\zeta, \zeta]$, where $\zeta$ is the noise magnitude and is set to $0.01$. Finally the vector $\y$ was obtained as $\y = A\x_* + \varepsilon$. 

We use the Matlab  package of NESTA~\footnote{\url{http://statweb.stanford.edu/~candes/nesta/}}. For fair comparison, we also use the projection code in the NESTA package for conducting projection. 
To handle the unknown smoothness parameter in the proposed algorithm, we use the backtracking technique~\citep{Beck:2009:FIS:1658360.1658364}. The parameter $\gamma$ is initially set to $0.001$ and decreased by half every $5000$ iterations after a  projection and the target smoothing parameter in NESTA is set to $10^{-5}$.  For the value of $\lambda$ in LoPNAG, we tune it from its theoretical value to several smaller values and choose the one that yields the fastest convergence.   We report the results  in Table~\ref{tab:cs}, which include different number of iterations, the corresponding number of projections, the recovery error of the found solution compared to the underlying true sparse solution, the objective value (i.e., the $\ell_1$ norm of the found solution) and the running time. Note that each iteration of NESTA requires two projections because it maintains two extra sequence of solutions.   From the results, we can see that LoPNAG converges significantly faster than NESTA. Even with only one projection, we are able to obtain a better solution than that of NESTA after running $10000$ iterations. 

\subsection{High-dimensional Distance Metric Learning}
\vspace*{-0.05in}
Consider the following distance metric learning problem: 
\begin{equation}\label{eqn:dml}
\min_{A\succeq 0}\frac{1}{2|\mathcal E|}\sum_{(i, j)\in\mathcal E}(1-y_{ij}- \|\x_i - \x_j\|_A^2)^2+ \tau\|A\|_{1}^{\text{off}},
\vspace*{-0.1in}
\end{equation}
where $\mathcal E$ denotes all pairs of training examples,  $y_{ij}=1$ indicates $\x_i, \x_j$ belong to the same class and $y_{ij}=-1$ indicates they belong to different classes, $\|\z\|_A^2 = \z^{\top}A\z$ and $\|A\|_{1}^{\text{off}}=\sum_{i\neq j}|A_{ij}|$. We note that such a formulation is useful for high dimensional problems due to the $\ell_1$ regularizer. A similar formulation with different forms of loss function has been adopted  in literature~\cite{QiTZCZ09}. We consider the square loss because it gives us faster convergence with a logarithmic number of projections by LoPGD.  
Due to the presence of the non-smooth PSD constraint and the $\ell_1$ regularizer, Nesterov's accelerated proximal gradient methods can not be applied efficiently to solving~(\ref{eqn:dml}) and the augmented unconstrained problem. Nevertheless, we can apply the proposed LoPGD method for solving the problem with a logarithmic number of projections. Regarding the constant $\theta$ in the local error bound condition for~\eqref{eqn:dml}, it still remains an open problem. Nonetheless,  a local error bound condition with $\theta=0.5$ might be established under certain regularity condition of the problem~\cite{ZhouSo15,CuiDing17}. For example, \citet{CuiDing17} provided a direct analysis of a local error bound condition with $\theta=0.5$ for a class of constrained convex symmetric matrix optimization problems regularized by nonsmooth spectral functions (including the indicator function of a PSD constraint). They established sufficient conditions (Theorem 16) for a local error bound condition with $\theta=0.5$ to hold, which  reduces to a regularity condition for~(\ref{eqn:dml}) depending on the optimal solutions of the problem. A thorough analysis of the regularity condition is much more involved and  left as an open problem. 

Next, we demonstrate the empirical performance of LoPGD for solving~(\ref{eqn:dml}).   
We use the colon-cancer data available on libsvm web portal, which  has 2000 features and  62 examples.   Fourty examples are used as training examples to generate 780 pairs to learn the distance metric. The regularization parameter is set to $\tau = 0.001$. We compare LoPGD, gradient descent method with only one projection (referred to as OPGD), and standard projected GD (referred to PGD). The step size in PGD and OPGD is set to $\eta_0/\sqrt{t}$, where $t$ is the iteration index.  We use the same tuned initial step size for all algorithms. The number of iterations per-epoch in LoPGD is set to 1000. The penalization parameter $\lambda$ in both OPGD and LoPGD is tuned and set to $10$.  In Table~\ref{tab:dml}, we report the objective values, the \#of iterations/projections, and running time across the first 8000 iterations. We can see that  LoPGD converges dramatically faster than PGD and also much faster than OPGD.

\vspace*{-0.05in}
\section{Conclusion}\label{sec:conc}
\vspace*{-0.05in}
 We have developed a general theory of optimization with only one projection for a family of  inequality constrained convex optimization problems. It yields an improved iteration complexity for smooth optimization compared with non-smooth optimization.  By exploring the local error bound condition, we further develop new algorithms with a logarithmic number of projections and achieve better convergence for both smooth and non-smooth optimization without strong convexity assumption. Applications in compressive sensing and distance metric learning demonstrate the effectiveness of the proposed improved algorithms. 

\section*{Acknowledgements}\label{sec:ack}
We are grateful to  all anonymous reviewers for their helpful comments. T. Yang is partially supported by National Science Foundation (IIS-1463988, IIS-1545995). L. Zhang thanks the support from NSFC (61603177) and JiangsuSF (BK20160658).

 \bibliography{all}
\bibliographystyle{icml2017}

\appendix
\section{Proof of Lemma~\ref{lem:0}}
When $c(\x)\leq 0$, $\x^\natural=\x$. There is nothing to prove. Therefore we consider $c(\x)>0$ and $\x^\natural\neq \x$. By KKT conditions, there exists $\zeta\geq 0$ and $\v\in\partial c(\x^\natural)$ such that 
\[
\x^\natural - \x + \zeta \v =0, \quad \text{ and }\quad \zeta c(\x^\natural)=0
\]
Since $\x^\natural\neq \x$, then $\zeta>0$, $c(\x^\natural)=0$ and $\v\neq 0$. 
 Therefore, $\x - \x^\natural $ is the same direction as $\v$. 
On the other hand, 
\begin{align*}
c(\x) &=c(\x) - c(\x^\natural) \geq (\x - \x^\natural)^{\top}\v\\
&= \|\v\| \|\x - \x^\natural\| \geq \rho \|\x - \x^\natural\|
\end{align*}
where the second equality uses the fact that $\x - \x^\natural$ is the same direction as $\v$ and the last inequality uses the inequality~(\ref{eqn:keyi0}).

\section{Proof of Proposition~\ref{prop:1}}
The two inequalities are straightforward to prove. We prove the smoothness property. Let $q(z)= \exp(z)/(1+\exp(z))$. It is not difficult to see that $q(z)$ is $1/4$-Lipschtiz continuous function. The gradient of $h_\gamma(\x)$ is given by  
\[
\nabla h_\gamma(\x) = \frac{\exp(\lambda c(\x)/\gamma)}{1+\exp(\lambda c(\x)/\gamma)}\lambda \nabla c(\x) = q(\lambda c(\x)/\gamma)\nabla c(\x)
 \]
 Then
 \begin{align*}
 &\left\|\nabla h_\gamma(\x)  - \nabla h_\gamma(\y) \right\| \\
 &=  \left\| q(\lambda c(\x)/\gamma)\lambda \nabla c(\x)  -   q(\lambda c(\y)/\gamma)\lambda \nabla c(\y)) \right\|\\
 &\leq  \left\| q(\lambda c(\x)/\gamma)\lambda \nabla c(\x)  -   q(\lambda c(\y)/\gamma)\lambda \nabla c(\x) \right\|\\
 & + \left\|  q(\lambda c(\y)/\gamma)\lambda \nabla c(\x) -   q(\lambda c(\y)/\gamma)\lambda \nabla c(\y)) \right\|\\
 &\leq \frac{\lambda G_C}{4}|\lambda c(\x)/\gamma - \lambda c(\y)/\gamma| + \lambda L_c\|\x - \y\|\\
 &\leq \left(\frac{\lambda^2G_c^2}{4\gamma}+\lambda L_c\right)\|\x - \y\|
 \end{align*}

\section{Proof of Corollary~\ref{cor:2}}
The following proposition shows the convergence of $F(\x)$. 
\begin{prop}~\citep{Beck:2009:FIS:1658360.1658364,opac-b1104789}
Assume $F(\x)$ is $L_F$-smooth. Let~(\ref{eqn:NAG}) run for $T$ iterations. If $F(\x)$ is convex, we can set $\beta_t = \frac{\tau_{t-1} -1}{\tau_t}$, where $\tau_t = \frac{1+\sqrt{1+4\tau^2_{t-1}}}{2}$ with $\tau_0=1$. Then for any $\x\in\R^d$ we  have 
\[
F(\x_T) - F(\x)\leq \frac{2L_F\|\y_0 - \x\|^2}{(T+1)^2}
\]
If $F(\x)$ is $\mu$-strongly convex, we can set $\beta_t = \frac{\sqrt{L_f} - \sqrt{\mu}}{\sqrt{L_f} + \sqrt{\mu}}$. Then for any $\x\in\R^d$ we have
\begin{align*}
&F(\x_T) - F(\x)\\
&\leq\exp\left(-T\sqrt{\frac{\mu}{L_F}}\right)\left( F(\y_0) - F(\x) + \frac{\mu}{2}\|\y_0 - \x\|^2\right)
\end{align*}
\end{prop}
We first prove the convergence for a smooth convex function $f(\x)$. From Theorem~\ref{thm:1} and the construction of $h_\gamma(\x)$ in~(\ref{eqn:pens}), we have
\[
f(\xt_T) - f(\x_*)\leq p \left(\gamma \ln 2 + \frac{2L_F\|\y_0 - \x_*\|^2}{(T+1)^2}\right)
\]
where $p = \frac{\lambda\rho}{\lambda\rho - G}$. 
By Proposition~\ref{prop:1}, we have $L_F =L_f +  \lambda L_c + \frac{\lambda^2G_c^2}{4\gamma}$. Then 
\begin{align*}
&f(\xt_T) - f(\x_*)\leq p \gamma \ln 2  \\
&+ p \left( \frac{\lambda^2 G_c^2\|\y_0 - \x_*\|^2}{2\gamma(T+1)^2} + \frac{2(L_f+\lambda L_c)\|\y_0 - \x_*\|^2}{(T+1)^2}\right)\\
&\leq  p\left(\gamma \ln 2 +  \frac{\lambda^2 G_c^2D^2}{2\gamma(T+1)^2} + \frac{2(L_f+\lambda L_c)D^2}{(T+1)^2}\right)\\
& = p\left( \frac{\sqrt{2\ln 2}\lambda  G_cD}{(T+1)} + \frac{2(L_f+\lambda L_c)D^2}{(T+1)^2}\right)
\end{align*}
where the last equality is due to the value of $\gamma$. 
Next, we prove the convergence for a smooth and strongly convex function $f(\x)$. First, we have
\begin{align*}
&F(\x_T) - F(\x_*)\leq \exp\left(-T\sqrt{\frac{\mu}{L_F}}\right)\cdot\\
&\left( \nabla F(\x_*)^{\top}(\y_0 - \x_*) + \frac{L_F+\mu}{2}\|\y_0 - \x_*\|^2\right)\\
&\leq \exp\left(-T\sqrt{\frac{\mu}{L_F}}\right)\left( \bar G\|\y_0 - \x_*\| + L_F\|\y_0 - \x_*\|^2\right)\\
&\leq \exp\left(-T\sqrt{\frac{\mu}{L_F}}\right) \bar GD +  \exp\left(-T\sqrt{\frac{\mu}{L_F}}\right)L_FD^2
\end{align*}
Note that $\x_*$ is not the optimal solution to $F(\x)$, hence $\nabla F(\x_*)\neq 0$ and we use its Lipschitz continuity property where $\bar G = G + \lambda G_c$. Following Theorem~\ref{thm:1}, we have 
\begin{align*}
&f(\xt_T) - f(\x_*)\leq p\gamma \ln 2 \\
&+p\left( \exp\left(-T\sqrt{\frac{\mu}{L_F}}\right) \bar GD +  \exp\left(-T\sqrt{\frac{\mu}{L_F}}\right)L_FD^2\right)\\
& \leq p \gamma \ln 2  \\
&+ p  \exp\left(-T\sqrt{\frac{\mu}{L_F}}\right) \bar GD+p \exp\left(-T\sqrt{\frac{\mu}{L_F}}\right)\frac{\lambda^2 G_c^2D^2}{4\gamma}\\
& +  p\left( \exp\left(-T\sqrt{\frac{\mu}{L_F}}\right)(L_f + \lambda L_c)D^2  \right)
\end{align*}
To avoid clutter, we will consider the dominating term. Consider $\gamma =\frac{1}{T^{2\alpha}}\leq 1$.   To bound the second term,  we have 
\begin{align*}
& \exp\left(-T\sqrt{\frac{\mu}{L_F}}\right)\frac{\lambda^2 G_c^2D^2}{4\gamma} = O\left( \exp\left(-T\sqrt{\frac{\mu}{L_F}}\right)T^{2\alpha}\right)\\
 & =  O\left( \exp\left(-T\sqrt{\frac{\mu}{L_f +\lambda L_c + \frac{\lambda^2G_c^2}{4\gamma}}}\right)T^{2\alpha}\right)\\
 & \leq  O\left(\exp\left(-T\sqrt{\frac{\mu \gamma}{L_f +\lambda L_c +\lambda^2G_c^2/4}}\right)T^{2\alpha}\right)\\
 & = O\left(\exp\left(-T^{1-\alpha}\sqrt{\frac{\mu }{L_f +\lambda L_c +\lambda^2G_c^2/4}}\right)T^{2\alpha}\right)
\end{align*}
Consider $T$ to be sufficiently larger such that $T\geq \left(\frac{L_f + \lambda L_c + \lambda^2 G_c^2/4}{\mu}\right)^{\frac{1}{2(1-\alpha)}}(4\alpha\ln T)^{\frac{1}{1-\alpha}}$, then 
\[
\exp\left(-T^{1-\alpha}\sqrt{\frac{\mu }{L_f +\lambda L_c +\lambda^2G_c^2/4}}\right)\leq\frac{ 1}{T^{4\alpha}}
\]
Therefore 
\[
 \exp\left(-T\sqrt{\frac{\mu}{L_F}}\right)\frac{\lambda^2 G_c^2D^2}{4\gamma} \leq O\left(\frac{1}{T^{2\alpha}}\right)
\]
Similarly we can show  the last two terms are dominated by $O(1/T^{4\alpha})$. 
As a result, 
\[
f(\xt_T) - f(\x_*)\leq O\left(\frac{1}{T^{2\alpha}} + \frac{1}{T^{4\alpha}}\right)
\]

\section{Proof of Theorem~\ref{thm:2}}
We first present a key lemma. 
\begin{lemma}\label{lem:key}
Let $D_{\epsilon}= \max_{\x\in\mathcal L_\epsilon}dist(\x,\Omega_*)$.  Then for any $\x\in\Omega$ and $\epsilon>0$ we have
\begin{align}\label{eqn:keyii}
\|\x - \x_\epsilon^\dagger\|\leq \frac{D_\epsilon}{\epsilon}(f(\x) - f(\x_\epsilon^\dagger)).
\end{align}
\end{lemma}
The above lemma was established  in~\citep{Yang15rsg}. 
\begin{proof}We consider $\x\not\in\S_\epsilon$, otherwise the conclusion holds trivially. 
By the first-order optimality conditions of~(\ref{eqn:ec}), we have for any $\u\in\Omega$, there exists $\zeta\geq 0$ (the Lagrangian multiplier of problem~(\ref{eqn:ec}))
\begin{equation}\label{eqn:o2}
\begin{aligned}
&(\x^\dagger_\epsilon - \x + \zeta \partial f(\x^\dagger_\epsilon))^{\top}(\u - \x_\epsilon^\dagger)\geq 0\\
\end{aligned}
\end{equation}
Let $\u = \x$ in the first inequality we have
\[
 \zeta \partial f(\x^\dagger_\epsilon)^{\top}(\x - \x^\dagger_\epsilon)\geq \|\x - \x_\epsilon^\dagger\|^2
\]
We argue that $\zeta>0$, otherwise $\x = \x_\epsilon^\dagger$ contradicting to the assumption $\x\not\in\S_\epsilon$. 
Therefore
\begin{align}\label{eqn:b1}
&f(\x) - f(\x^\dagger_\epsilon)\geq  \partial f(\x^\dagger_\epsilon)^{\top}(\x - \x^\dagger_\epsilon)\geq \frac{\|\x - \x_\epsilon^\dagger\|^2}{\zeta}\nonumber\\
&=  \frac{\|\x - \x_\epsilon^\dagger\|}{\zeta}\|\x - \x_\epsilon^\dagger\|
\end{align}
Next we prove that $\zeta$ is upper bounded. Since
\begin{align*}
-\epsilon = f(\x^*_\epsilon) - f(\x^\dagger_\epsilon)\geq (\x_\epsilon^* - \x^\dagger_\epsilon)^{\top}\partial f(\x_\epsilon^\dagger)
\end{align*}
where $\x^*_\epsilon$ is the closest point to $\x^\dagger_\epsilon$ in the optimal set. Let $\u=\x^*_\epsilon$ in the inequality of~(\ref{eqn:o2}), we have
\begin{align*}
(\x_\epsilon^\dagger - \x)^{\top}(\x^*_\epsilon - \x^\dagger_\epsilon)\geq \zeta ( \x^\dagger_\epsilon - \x^*_\epsilon)^{\top}\partial f(\x^\dagger_\epsilon)\geq \zeta \epsilon
\end{align*}
Thus 
\[
\zeta \leq \frac{(\x_\epsilon^\dagger - \x)^{\top}(\x^*_\epsilon - \x^\dagger_\epsilon)}{\epsilon}\leq \frac{D_{\epsilon}\|\x_\epsilon^\dagger - \x\|}{\epsilon}
\]
Therefore 
\[
\frac{\|\x - \x_\epsilon^\dagger\|}{\zeta}\geq \frac{\epsilon}{D_\epsilon}
\]
Combining the above inequality with~(\ref{eqn:b1}) we have
\[
f(\x) - f(\x^\dagger_\epsilon)\geq \frac{\epsilon}{D_\epsilon}\|\x - \x_\epsilon^\dagger\|
\]
which completes the proof. 
\end{proof}

\begin{prop}[\citet{zinkevich-2003-online}]\label{prop:gd}
Let $\x_{t+1} = \x_t - \eta\partial F(\x_t)$ run for $T$ iterations. Assume $\|\partial F(\x)\|\leq \bar G$. Then for any $\x\in\R^d$
\begin{align*}
F(\xh_T) - F(\x)\leq  \frac{\eta \bar G^2}{2} + \frac{\|\x_1 - \x\|^2}{2\eta T}
\end{align*}
where $\xh_T = \sum_{t=1}^T\x_t/T$. 
\end{prop}

\subsection{Proof of Theorem~\ref{thm:2}}
Let $\epsilon_k = \frac{\epsilon_0}{2^k}$. We assume $\x_1,\ldots, \x_{K-1}\not\in\S_{2\epsilon}$; otherwise the result holds trivially. Let $\x_{k,\epsilon}^\dagger\in\Omega$ denote the closest point to $\x_k$ in the sublevel set $\S_\epsilon$ of $f(\x)$. Then $f(\x_{k,\epsilon}^\dagger)=f_* + \epsilon, k=1,\ldots, K-1$, which is because we assume that  $\x_1,\ldots, \x_{K-1}\not\in\S_{2\epsilon}$. 
We will prove by induction that $f(\x_k) - f_* \leq \epsilon_k + \epsilon$. It holds for $k=0$ because of Assumption (i). We assume it holds for $k-1$ and prove it is true for $k\leq K$. 

We consider the $k$-th epoch of LoPGD. Applying Proposition~\ref{prop:gd}, we have for any $\x\in\R^d$
\[
F(\xh_k) - F(\x)\leq \frac{\eta_k \bar G^2}{2} + \frac{\|\x_{k-1} - \x\|^2}{2\eta_k t}
\]
 Let $\x = \x_{k-1,\epsilon}^\dagger\in\Omega$. Then 
\begin{align*}
F(\xh_k) - F(\x_{k-1,\epsilon}^\dagger)\leq  \frac{\eta_k \bar G^2}{2} + \frac{\|\x_{k-1} - \x_{k-1,\epsilon}^\dagger\|^2}{2\eta_k t}
\end{align*}
Since 
\begin{align*}
F(\x_{k-1,\epsilon}^\dagger)& = f(\x_{k-1,\epsilon}^\dagger) + \lambda [c(\x_{k-1,\epsilon}^\dagger)]_+ =f(\x_{k-1,\epsilon}^\dagger) \\
F(\xh_k)&= f(\xh_k) + \lambda[c(\xh_k)]_+
\end{align*}
Then
\begin{align*}
&f(\xh_k) + \lambda[c(\xh_k)]_+ - f(\x_{k-1,\epsilon}^\dagger)\\
&\leq  \underbrace{ \left(\frac{\eta_k \bar G^2}{2} + \frac{\|\x_{k-1} - \x_{k-1,\epsilon}^\dagger\|^2}{2\eta_k t}\right)}\limits_{B_t}
\end{align*}
Then 
\begin{align*}
f(\xh_k) + \lambda [c(\xh_k)]_+\le f(\x_{k-1,\epsilon}^\dagger) + B_t
\end{align*}
Then 
\begin{align*}
&\lambda \rho \|\xh_k - \x_k\|\leq f(\x_{k-1,\epsilon}^\dagger) - f(\xh_k) + B_t\\
&\leq f(\x_{k-1,\epsilon}^\dagger) - f(\x_k) + f(\x_k)  - f(\xh_k) + B_t
\end{align*}
Assume $f(\x_k) - f_*> \epsilon$ (otherwise the proof is done), thus $f(\x_{k-1,\epsilon}^\dagger)\leq f(\x_k)$. Then
\begin{align*}
\lambda \rho \|\xh_k - \x_k\|&\le G\|\x_k  -\xh_k\| + B_t
\end{align*}
leading to 
\begin{align*}
\|\xh_k - \x_k\|\leq \frac{B_t}{\lambda\rho - G}
\end{align*}
Then
\begin{align*}
f(\x_k) &- f(\x_{k-1,\epsilon}^\dagger)\le f(\x_k) - f(\xh_k) + f(\xh_k) - f(\x_{k-1,\epsilon}^\dagger)\\
&\le G\|\xh_k - \x_k\| + B_t = \frac{\lambda\rho}{\lambda\rho - G}B_t \\
&= p \left(\frac{\eta_k \bar G^2}{2} + \frac{\|\x_{k-1} - \x_{k-1,\epsilon}^\dagger\|^2}{2\eta_k t}\right)\\
&\leq p \left(\frac{\eta_k \bar G^2}{2} + \frac{D_\epsilon^2(f(\x_{k-1}) - f(\x_{k-1,\epsilon}^\dagger))^2}{2\eta_k t\epsilon^2}\right)\\
&\leq  p \left(\frac{\eta_k \bar G^2}{2} + \frac{\sigma^2(f(\x_{k-1}) - f(\x_{k-1,\epsilon}^\dagger))^2}{2\eta_k t\epsilon^{2(1-\theta))}}\right)
\end{align*}
where the third inequality uses Lemma~\ref{lem:key} and the last inequality uses the local error bound condition. 
Since we assume $f(\x_{k-1}) - f_*\leq \epsilon_{k-1} + \epsilon$. Thus $f(\x_{k-1}) - f(\x_{k-1,\epsilon}^\dagger) \leq \epsilon_{k-1}$.  Then
\begin{align*}
f(\x_k) - f(\x_{k-1,\epsilon}^\dagger)&  \leq \frac{p\eta_k \bar G^2}{2} + \frac{p\sigma^2\epsilon_{k-1}^2}{2\eta_k t\epsilon^{2(1-\theta))}}
\end{align*}
By noting the values of  $\eta_k = \frac{\epsilon_{k-1}}{2p\bar G^2}$ and $t = \frac{4\sigma^2p^2\bar G^2}{\epsilon^{2(1-\theta)}}$, we have
\begin{align*}
f(\x_k) - f(\x_{k-1,\epsilon}^\dagger)&\leq \frac{\epsilon_{k-1}}{4} + \frac{\epsilon_{k-1}}{4} = \epsilon_k. 
\end{align*} 
Therefore
\begin{align*}
f(\x_k) - f_* \leq \epsilon_k + \epsilon
\end{align*}
due to the assumption $f(\x_{k-1})\geq f_* + 2\epsilon$ and $f(\x_{k-1,\epsilon}^\dagger)=\epsilon$. By induction, we therefore show that with at most $K=\log_2(\epsilon_0/\epsilon)$ epochs, we have
\[
f(\x_K) - f_*\leq \epsilon_K +\epsilon \leq 2\epsilon
\]

\section{Proof of Theorem~\ref{thm:3}}
Following a similar  analysis and using the convergence of NAG and  Proposition~\ref{prop:1}, we have
\begin{align*}
f(\x_k) - f(\x_{k-1,\epsilon}^\dagger)&\le p \gamma_k\ln 2\\
&+p\left(  \frac{\lambda^2 G_c^2\|\x_{k-1} - \x_{k-1,\epsilon}^\dagger\|^2}{2\gamma_kt_k^2}\right)\\
&+p\frac{2(L_f+\lambda L_c) \|\x_{k-1} - \x_{k-1,\epsilon}^\dagger\|^2}{t_k^2}
\end{align*}
By using Lemma~\ref{lem:key} and the local error bound condition, we have
\begin{align*}
f(\x_k) - f(\x_{k-1,\epsilon}^\dagger)&\le p \left(\gamma_k\ln 2 + \frac{\lambda^2 G_c^2\sigma^2\epsilon_{k-1}^2}{2\gamma_kt_k^2\epsilon^{2(1-\theta)}}\right)\\
&+ p  \frac{2(L_f+\lambda L_c)\sigma^2\epsilon_{k-1}^2}{t_k^2\epsilon^{2(1-\theta)}}
\end{align*}
Plugging  the values of  $\gamma_k  = \frac{\epsilon_{k-1}}{6p\ln 2 }$, $t_k =\frac{\sigma}{\epsilon^{1-\theta}} \max \{\lambda G_cp\sqrt{18\ln 2}, \sqrt{12(L_f + \lambda L_c)\epsilon_{k-1}}\}$ into the above inequality yields  
\[
f(\x_k) - f(\x_{k-1,\epsilon}^\dagger)\leq \frac{\epsilon_{k-1}}{6} + \frac{\epsilon_{k-1}}{6} + \frac{\epsilon_{k-1}}{6} = \epsilon_k
\]
Then 
\[
f(\x_k) - f_*\leq \epsilon_k + \epsilon
\]
Therefore the total number of iterations is 
{\footnotesize\begin{align*}
&\frac{\sigma}{\epsilon^{1-\theta}}\max\left( \lambda G_cp\sqrt{18\ln 2}\log_{2}(\epsilon_0/\epsilon), \sqrt{12\bar L}\sum_{k=1}^K\epsilon^{1/2}_{k-1} \right)\\
& \leq\frac{\sigma}{\epsilon^{1-\theta}} \max\left( \lambda G_cp\sqrt{18\ln 2}\log_{2}(\epsilon_0/\epsilon), \sqrt{12\bar L}\frac{\sqrt{2\epsilon_0}}{\sqrt{2} - 1} \right) 
\end{align*}}

\section{Lower bound of the subgradient of the constraint function for a PSD constraint}
We first show that
\begin{align*}
&\text{Conv}\{-\u\u^{\top}| \|\u\|= 1, X\u = 0\}= \\
&   \text{Conv}\{-U|U\succeq 0,\text{Tr}(X^{\top}U)=0, \text{rank}(U)= 1, \text{Tr}(U)=1\}.
    \end{align*}
    
In fact, given any $\u\in\R^d$ with $ \|\u\|= 1$ and $X\u = 0$, we can show  $\u\u^{\top}\succeq 0$, $\text{Tr}(X^{\top}\u\u^{\top})=\u^{\top}X\u=0$, $\text{rank}(\u\u^{\top})= 1$ and $\text{Tr}(\u\u^{\top})=\|\u\|^2=1$, which means $\u\u^{\top}$ belongs to the set on the right. Since the set on the left is the convex hull of all such $\u\u^{\top}$, the set on the left is included in the set of the right.

On the other hand, given any element $U$ from the set of the right, we can represent it as 
$U=\sum_{k=1}^K\lambda_kU_k$ where $\sum_{k=1}^K\lambda_k=1$, $\lambda_k\geq0$, $U_k=\u_k\u_k^{\top}$ for some $\u_k$, $\text{Tr}(X^{\top}U_k)=0$ and $\text{Tr}(U_k)=1$ for $k=1,\dots,K$. These three properties of $U_k$ imply 
$X\u_k = 0$ and $\|\u_k\|^2=1$ so that $U$ is a convex combination of some elements of the set on the left. Therefore, the set on the right is included in the set of the left.

Next, we want to show 
    \begin{align*}
&\text{Conv}\{-U|U\succeq 0,\text{Tr}(X^{\top}U)=0, \text{rank}(U)= 1, \text{Tr}(U)=1\}
    \\
    &=\{-U|U\succeq 0,\text{Tr}(X^{\top}U)=0, \text{Tr}(U)=1\}
    \end{align*}
It is easy to see that the set on the left is a subset of  the set on the right. To show the opposite, given any element $U$ from the set of the right, we consider its eigenvalue decomposition 
$U=\sum_{k=1}^K\lambda_k\u_k\u_k^{\top}$ where $K\leq d$ and $\lambda_k>0$ and $\u_k$ are the eigenvalue and the corresponding eigenvector with $\|\u_k\|=1$. Since $X$ is PSD, the property $\text{Tr}(X^{\top}U)=0$ implies $\sum_{k=1}^K\lambda_k\u_k^{\top}X\u_k=0$ so that $X\u_k$ must be zero for $k=1,\dots,K$. As a result, $U=\sum_{k=1}^K\lambda_kU_k$ with $U_k=\u_k\u_k^{\top}$ being an element in the set on the left. Note that $\text{Tr}(U)=\sum_{k=1}^K\lambda_k\u_k^{\top}\u_k=\sum_{k=1}^K\lambda_k=1$. This means $U$ is in the set on the left also.

If the dimension of the null space of $X$ is $r$ with $1\leq r\leq d$ then we can write $X=V\Sigma V^{\top}$, where $\Sigma=\left(\begin{array}{cc}\Sigma_r&0\\ 0 & 0\end{array}\right)$ is a diagonal matrix with $\Sigma_r\in\R^{d- r, d- r}$. We can set the constant $\rho$ to be the solution of the following optimization problem. 
    \begin{align*}
    \rho =&\quad\arg\min_{U\in\mathbb{R}^{d\times d}} \|U\|_F\\
    s.t.&\quad U\succeq 0,\text{Tr}(X^{\top}U)=0, \text{Tr}(U)=1
    \end{align*}
To simplify the problem, we note that \[
\text{Tr}(X^TU)=\text{Tr}(V\Sigma V^{\top}U)=\text{Tr}(\Sigma V^{\top}UV)=0.\]
Let  $V^{\top}UV=\left(\begin{array}{cc}U_{11},&U_{12}\\ U_{21}, & U_{22}\end{array}\right)$ where $U_{11}\in\R^{(d-r)\times (d-r)}$ and $U_{22}\in\R^{r\times r}$. Because $\Sigma$ is a diagonal matrix with nonnegative entries, it then leads to that the diagonal entries of $U_{11}$ are all zeros, as a result  $U_{11}=0$ and consequentially $U_{21}=U_{12}=0$ due to that $V^{\top}UV\succeq 0$. As a result, $\|U\|_F = \|V^{\top}UV\|_F = \|U_{22}\|_F$ and $\text{Tr}(U)=\text{Tr}(V^{\top}UV) =\text{Tr}(U_{22})$. Therefore, we get
\begin{align*}
    \rho &=\quad\arg\min_{U_{22}\in\mathbb{R}^{r\times r}} \|U\|_F,\quad s.t.\quad U_{22}\succeq 0,\: \text{Tr}(U_{22})=1
\end{align*}
As a result, $\rho= \frac{1}{\sqrt{r}}\geq \frac{1}{\sqrt{d}}$.

\section{Proof of Theorem~\ref{thm:cs}}
The  proof follows similarly to that of Theorem 3.5 in~\citep{guoyincalculus2016} but tailored to the our problem at hand. Let $\x_*$ denote an optimal solution to $\min_{c(\x)\leq 0}f(\x)$. Let $\Omega=\{\x: c(\x)\leq 0\}$ and $F(\x) = f(\x) + \mathbb I_{\Omega}(\x)$, where $\mathbb I_{\Omega}(\x)$ denotes the indicator function associated with $\Omega$. 

Due to the strict feasibility condition, i.e., there exists $\x_0$ such that $\|A\x_0 - \y\|^2<\tau$, then by the Lagrangian theory there exists $\zeta\geq 0$ such that 
\[
F(\x_*) = \min_{\x\in\R^d}F(\x) = \min_{\x\in\R^d} f(\x) + \zeta c(\x) = f(\x_*) + \zeta c(\x_*)
\]
The KKT condition implies $0\in\partial f(\x_*) + \zeta \nabla c(\x_*)$ and $\zeta c(\x_*)=0$. Since $\x_*\neq 0$, thus $0\not\in\partial f(\x_*)$, and consequentially $\zeta>0$. As a result, $c(\x_*)=0$. In view of~\citep{rockafellar1970convex}[Theorem 28.1], we have
\[
\x_* \in\arg\min F(\x) =  \{\x: c(\x)=0\}\cap\arg\min f(\x) + \zeta c(\x)
\]
By noting the form of $f(\x)  + \zeta c(\x) = \|\x\|_1 + \zeta(\|A\x - \y\|^2 - \tau) = h(A\x) + \|\x\|_1$,  where $h(\cdot)$ is a strongly convex function, it must hold that $A\x$ is a constant for $\x\in\arg\min f(\x) + \zeta c(\x)$, i.e., $c(\x)$ is a constant over $\arg\min_{\x}f(\x) + \zeta c(\x)$. Then for any $\x\in\arg\min f(\x) + \zeta c(\x)$, we have $c(\x) = c(\x_*) = 0$. It follows that 
\begin{align}\label{eqn:objF}
&\arg\min F(\x) = \{\x: c(\x)=0\}\cap\arg\min f(\x) + \zeta c(\x)\notag\\
& = \arg\min f(\x) + \zeta c(\x)
\end{align}
Next, we consider the problem $\min_{\x}f_\zeta(\x) = f(\x) + \zeta c(\x)= \|\x\|_1 + \zeta(\|A\x - \y\|^2 - \tau)$. It has been shown that there exists $\sigma>0$ (e.g.~\citep{DBLP:journals/corr/nesterov16linearnon}[Theorem 4.3]) such that 
\begin{align}\label{eqn:eb}
\text{dist}(\u, \arg\min f_\zeta(\u))\leq \sigma (f_\zeta(\u) - \min_{\u}f_\zeta(\u))^{1/2}
\end{align}
for any $\u\in\R^d$ such that $f_\zeta(\u)-\min_{\u}f_\zeta(\u) \leq \epsilon$. Next, consider any $\x$ such that $c(\x)\leq 0$ and $f(\x)\leq f_* + \epsilon$. We have $f_\zeta(\x) - \min_\u f_\zeta(\u) = f(\x) + \zeta c(\x) - f(\x_*)\leq f(\x) - f_*\leq \epsilon$. As a result, (\ref{eqn:eb}) holds for any $\x$ such that $c(\x)\leq 0$ and $f(\x)\leq f_* + \epsilon$, i.e., 
\begin{align*}
&\text{dist}(\x, \arg\min f_\zeta(\u))\leq \sigma (f_\zeta(\x) - \min_{\u}f_\zeta(\u))^{1/2}\\
& = \sigma (f(\x) + \zeta c(\x) - f_*)^{1/2} \leq  \sigma (f(\x) - f_*)^{1/2}
\end{align*}
In view of~(\ref{eqn:objF}), we can finish the proof. 

\end{document}